\documentclass{amsart}
\date{\today}
\usepackage{amssymb}
\usepackage{latexsym}
\newcommand{\Z}{{\mathbb Z}}
\newcommand{\R}{{\mathbb R}}
\newcommand{\C}{{\mathbb C}}

\newcommand{\D}{{\mathbb D}}

\def\H{{\mathbb {H}}}

\def\SL{\mathrm{SL}}

\def\dist{\operatorname{dist}}

\def\supp{\operatorname{supp}}

\newtheorem{theorem}{Theorem}
\newtheorem{lemma}{Lemma}

\theoremstyle{remark}
\newtheorem{rem}{Remark}[section]

\def\bm{\begin{pmatrix}}
\def\em{\end{pmatrix}}

\begin{document}

\title{Generic Singular Spectrum For Ergodic Schr\"odinger Operators}

\author{Artur Avila and David Damanik}

\address{
Laboratoire de Probabilit\'es et Mod\`eles al\'eatoires\\
Universit\'e Pierre et Marie Curie--Boite courrier 188\\
75252--Paris Cedex 05, France
}
\email{artur@ccr.jussieu.fr}

\address{
Mathematics 253-37\\
California Institute of Technology\\
Pasadena, CA 91125, USA }
\email{damanik@caltech.edu}

\thanks{D.\ D.\ was supported in part by NSF grant DMS--0227289.}

\begin{abstract}

We consider Schr\"odinger operators with ergodic potential $V_\omega(n)=f(T^n(\omega))$,
$n \in \Z$, $\omega \in \Omega$, where $T:\Omega \to \Omega$ is a non-periodic
homeomorphism.  We show that for generic $f \in C(\Omega)$, the spectrum has no
absolutely continuous component. The proof is based on approximation by discontinuous
potentials which can be treated via Kotani Theory.

\end{abstract}

\maketitle

\section{Introduction}

Let $\Omega$ be a compact metric space, $T:\Omega \to \Omega$ a homeomorphism, and $\mu$
a $T$-ergodic Borel measure.  We will always assume that $T$ is not periodic, that is,
$\mu$ is non-atomic.  For a bounded and measurable function $f : \Omega \to \R$, we
consider (line) Schr\"odinger operators $H_\omega = \Delta + V_\omega$, $n \in \Z$ with
potential $V_\omega(n) = f( T^n \omega )$ and the associated Lyapunov exponents $\gamma
(z)$, $z \in \C$. By Kunz-Souillard (cf.\ \cite{cfks,ks}), there exists a compact set
$\Sigma_{{\rm ac}}(f) \subseteq \R$ such that $\sigma_{{\rm ac}}(H_\omega) = \Sigma_{{\rm
ac}}(f)$ for $\mu$-a.e.\ $\omega \in \Omega$. By Pastur-Ishii-Kotani (cf.\
\cite{cfks,i,k,p,s}), $\Sigma_{{\rm ac}}(f) = \overline{\{ E \in \R : \gamma(E) = 0
\}}^{{\rm ess}}$.

We shall only consider situations where the potentials $V_\omega$ are not periodic. In
this case, it is an interesting question whether there can be any absolutely continuous
spectrum.

It was observed by Kotani, \cite{k1}, that $\Sigma_{{\rm ac}}(f)$ is empty if $f$ takes
only finitely many values. Damanik and Killip, \cite{dk}, derived the same conclusion
under the assumption that $f$ is discontinuous at some point $\omega_0$ (but continuous
at all points in the forward orbit of $\omega_0$ under $T$). Here, we will consider the
case of continuous $f$. It is not in general true in this case that $\Sigma_{{\rm
ac}}(f)$ is empty: If $\Omega = \R / \Z$, $T$ is a Diophantine irrational rotation, $v$
is analytic and $\lambda > 0$ is small enough, then for $f = \lambda v$, the spectrum is
(almost surely) purely absolutely continuous; see, for example, Bourgain and Jitomirskaya
\cite{bj} and references therein. However, we show that absence of absolutely continuous
spectrum is a generic phenomenon:

\begin{theorem}\label{main}
There is a residual set of functions $f$ in $C(\Omega)$ such that $\Sigma_{{\rm ac}}(f) =
\emptyset$.
\end{theorem}

\begin{rem}

A subset of a complete metric space (or more generally a Baire space)
is called residual or (Baire) generic if it contains a countable
intersection of dense open sets. By Baire's Theorem, such a set is dense.

\end{rem}

\begin{rem}

It is well known (see \cite {B} for a more general result) that for generic potentials,
the Lyapunov exponent vanishes generically in the spectrum.  Frequently (for instance, if
$T:\Omega \to \Omega$ has a periodic point in the support of $\mu$), the spectrum
contains intervals in a persistent way, and hence has positive Lebesgue measure.  This is
not incompatible with our result, since the Lyapunov exponent can be discontinuous.

\end{rem}

The method used in the proof of Theorem~\ref {main} can be applied to certain
one-parameter families of potentials.  This has the following consequence.

\begin{theorem}\label{main2}
There is a residual set of functions $f$ in $C(\Omega)$ such that $\Sigma_{{\rm
ac}}(\lambda f) = \emptyset$ for almost every $\lambda > 0$.
\end{theorem}

This is particularly striking in the case of quasiperiodic systems.  In this context one
often expects absolutely continuous spectrum for small coupling, and as we mentioned
before, this has been established under strong regularity assumptions on the potential
(Bourgain-Jitomirkaya's result).  It was not clear whether the assumption on the
regularity was an artifact of the known methods.  Our result shows that some regularity
is indeed necessary: it is not enough to assume only continuity of the potential.  It is
an interesting open problem to determine the weakest possible regularity assumption on
$f$ under which a Bourgain-Jitomirskaya-type result holds.

\medskip

\noindent\textit{Acknowledgements.} This work was done while the first author was
visiting Caltech. We would like to thank Svetlana Jitomirskaya and Barry Simon for
stimulating discussions.

\section{A Semi-Continuity Result}

We will need some basic facts about Lyapunov exponents and $m$-functions,
see \cite {s}.  The Lyapunov exponent is defined by
$$
\gamma_f(E)=\lim_{n \to \infty} \frac {1} {n} \int_\Omega \ln \|S_{f,E}^n(\omega)\| \,
d\mu(\omega),
$$
where
$$
S_{f,E}^n(\omega)=S_{f,E}(T^{n-1}(\omega)) \cdots S_{f,E}(\omega),
$$
and
$$
S_{f,E}(\omega)=\bm E-f(\omega) & -1\\1 & 0 \em \in \SL(2,\C).
$$
We have that $E \mapsto \gamma_f$ is a non-negative real-symmetric subharmonic function.
If $E \in \H=\{z \in \C : \Im z>0\}$, we have the formula
$$
\gamma_f(E)=\int_\Omega -\Re \ln m_{\omega,f}(E) \, d\mu(\omega),
$$
where $m_{\omega,f}:\H \to \H$ is a holomorphic function, called
the $m$-function, which is given almost everywhere by
$$
m_{\omega,f}(E)=\lim_{n \to \infty} S_{f,E}^n(T^{-n}(\omega)) \cdot i
$$
(here we consider the usual action of $\SL(2,\C)$ by M\"obius transformations, i.e., $\bm
a&b\\c&d \em \cdot z=\frac {az+b} {cz+d}$).

Define $M(f) = | \{ E \in \R : \gamma_f(E) = 0 \} |$, where $| \cdot |$ denotes Lebesgue
measure. As a consequence of the results of Ishii, Kotani, and Pastur, we have
$\Sigma_{{\rm ac}}(f) = \emptyset$ if and only if $M(f) = 0$.

\begin{lemma}\label{L1}
For every $r > 0$, $\Lambda>0$, the maps
\begin{equation} \label {L11}
(L^1(\Omega) \cap B_r(L^\infty(\Omega)), \| \cdot \|_1) \to \R, \;
\; f \mapsto M(f)
\end{equation}
and
\begin{equation} \label {L12}
(L^1(\Omega) \cap B_r(L^\infty(\Omega)), \| \cdot \|_1) \to \R, \;
\; f \mapsto \int_0^\Lambda M(\lambda f) d\lambda
\end{equation}
are upper semi-continuous.

\end{lemma}

\begin{proof}

It is enough to show that (\ref {L11}) is upper semi-continuous, since this
implies that (\ref {L12}) is also upper semi-continuous by Fatou's Lemma.

We have to show that if $(f_n)_{n \in \Z_+}$, $f$ are uniformly
bounded in $L^\infty$ and $f_n \to f$ in $L^1$, then $\limsup
M(f_n) \le M(f)$.

Assume otherwise. Then (by passing to a suitable subsequence),
there are a constant $C < \infty$ and a sequence $(f_n)$ such that
\begin{itemize}

\item[(i)] $f_n \to f$ in $L^1$ and pointwise,

\item[(ii)] $\|f_n\|_\infty \le C$, $\|f\|_\infty \le C$,

\item[(iii)] $\liminf M(f_n) \ge M(f) + \varepsilon$ for some
$\varepsilon > 0$.

\end{itemize}

By (i), we have pointwise convergence of the m-functions
$m_{\omega,f_n}$ in $\H$ for almost every $\omega$. Thus, by
dominated convergence and (ii),
the associated Lyapunov exponents $\gamma_{f_n}(E)$
converge pointwise in $\H$ to $\gamma_f(E)$.

By (ii), all Lyapunov exponents are positive outside the interval $I = [-2-C,2+C]$. Thus,
we can limit our attention to this interval. Consider the region $U$ in $\H$, bounded by
the equilateral triangle $T$ with sides $I,J,K$. Consider a conformal mapping $\Phi$ from
the unit disk $\D$ to $U$. By the Schwarz-Christoffel formula (see, e.g., \cite{dt}),
\begin{equation}\label{SC}
\Phi'(z) = {\rm const} \cdot \prod_{j=1}^3 \left( 1 -
\tfrac{z}{z_j} \right)^{-2/3},
\end{equation}
where $z_1,z_2,z_3$ are the inverse images under $\Phi$ of the
vertices of $T$.

The functions $\gamma_{f_n} \circ \Phi$ are harmonic and bounded in
$\D$. This yields
$$
\gamma_{f_n}( \Phi(0) ) = \frac{1}{2 \pi} \int_0^{2\pi} \gamma_{f_n}
\left(\Phi(e^{i \theta})\right) \, d\theta,
$$
and similarly for $\gamma_f$. Since $\gamma_{f_n}( \Phi(0) ) \to
\gamma_f(\Phi(0) )$ as $n \to \infty$, we infer
$$
\frac{1}{2 \pi} \int_0^{2\pi} \left[ \gamma_{f_n} \left(\Phi(e^{i
\theta})\right) - \gamma_f \left(\Phi(e^{i \theta})\right) \right]
\, d\theta \to 0.
$$
By dominated convergence, the integrals along $J$ and $K$ go to
zero individually. Therefore,
$$
\int_I [\gamma_{f_n}(E) - \gamma_f(E)] g(E) \, dE \to 0
$$
where $g(E) = [\Phi'(\Phi^{-1}(E))]^{-1}$. It follows from
\eqref{SC} that $g$ vanishes at the endpoints of $I$ and is
continuous and non-vanishing inside $I$.

By upper semi-continuity of the Lyapunov exponent and dominated
convergence,
$$
\int_I \max \{\gamma_{f_n}(E) - \gamma_f(E),0\} \, g(E) \, dE \to
0,
$$
and hence
$$
\int_I \min\{\gamma_{f_n}(E) - \gamma_f(E),0\} \, g(E) \, dE \to
0.
$$
Consequently, since $\gamma_f|I$ is bounded and $\gamma_{f_n}|I$ is
non-negative,
$$
\int_I \min \{\gamma_{f_n}(E) - \gamma_f(E),0\} \, dE \to 0.
$$

Choose $\delta > 0$ such that the set $X = \{ E \in I : \gamma_f(E)
< \delta \}$ has measure bounded by $M(f) +
\frac{\varepsilon}{4}$, with $\varepsilon$ from (iii). Then,
$$
\int_{I \setminus X} \min\{\gamma_{f_n}(E) - \gamma_f(E),0\} \,
dE \to 0.
$$
This shows that for $n \ge n_0$, there exists a set $Y_n$ of
measure bounded by $\frac{\varepsilon}{4}$ such that $\gamma_{f_n}(E)
\ge \frac{\delta}{2}$ for every $E \in I \setminus (X \cup Y_n)$.
Consequently, $\limsup M(f_n) \le M(f) + \frac{\varepsilon}{2}$,
which contradicts (iii).
\end{proof}

\section{Approximation by Discontinuous Potentials}

\begin{lemma} \label {L3}
There exists a dense subset $\mathcal{Z}$ of $L^\infty(\Omega)$ such that if $s \in Z$,
then
\begin{enumerate}

\item $s(\omega)$, $\omega \in \Omega$, takes finitely many values,

\item $s(T^n(\omega))$, $n \in \Z$, is not periodic for almost every $\omega \in
\Omega$.

\end{enumerate}

\end{lemma}

\begin{proof}
Let $W_k$ be the closed subspace of functions $s$ taking at most $k$ values. Obviously
$W=\cup_{k \geq 2} W_k$ is dense in $L^\infty(\Omega)$.  So we only have to show that
there is a dense subset $S_k \subset W_k$ of functions satisfying the second property.
Given $s \in W$, $\omega \in \Omega$, let $\phi(s,\omega) \in \Z_+ \cup \{\infty\}$ be
the period of $s(T^n(\omega))$, $n \in \Z$.  Then $\phi(s,\omega)$ is a constant
$\Phi(s)$ almost everywhere.  Let $W_{k,n}=\{s \in W_k : \Phi(s) \leq n\}$.  It is easy
to see that $W_{k,n}$ is a closed subset of $W_k$ and $W_k \neq W_{k,n}$.  Thus, $W_k
\setminus \cup_{n \in \Z_+} W_{k,n}$ is dense in $W_k$.
\end{proof}

\begin{lemma}\label{L2}
For $f \in C(\Omega)$, $\varepsilon > 0$, $\delta>0$, $\Lambda>0$,
there exists $\tilde{f} \in C(\Omega)$ such that
$\| f - \tilde{f} \|_\infty < \varepsilon$, $M(\tilde{f}) < \delta$, and
$\int_0^\Lambda M(\lambda \tilde {f}) d\lambda<\delta$.
\end{lemma}

\begin{proof}
Let $\mathcal{Z}$ be as in Lemma \ref {L3} and choose $s \in \mathcal{Z}$ such that $\|f
- s\|_\infty < \frac{\varepsilon}{2}$.  By the Kotani result, \cite{k1}, we have
$M(\lambda s) = 0$ for every $\lambda>0$.  Next we choose continuous functions $f_n$, for
which we have $\|s - f_n\|_\infty < \frac{\varepsilon}{2}$ for all $n$ and $\|s - f_n\|_1
\to 0$ as $n \to \infty$.  For instance, take $f_n(\omega)=\int_\Omega C_n(\omega)^{-1}
c_n(\omega,\omega') s(\omega') d\nu(\omega')$, where $\nu$ is a probability measure with
$\supp \nu=\Omega$, $C_n(\omega)=\int_\Omega c_n(\omega,\omega') d\nu(\omega')$,
$c_n(\omega,\omega')=\max \{(n+n_0)^{-1}-\dist(\omega,\omega'),0\}$, $n_0$ sufficiently
large. Lemma~\ref{L1} implies $M(f_n), \int_0^\Lambda M(\lambda f_n) d\lambda \to 0$ as
$n \to \infty$. Thus, choosing $n$ large enough so that $M(f_n), \int_0^\Lambda M(\lambda
f_n) d\lambda<\delta$, we complete the proof.
\end{proof}

\begin{proof}[Proof of Theorem~\ref{main}.]
For $\delta > 0$, we define
$$
M_\delta = \{ f \in C(\Omega) : M(f) < \delta \}.
$$
By Lemma~\ref{L1}, $M_\delta$ is open, and by Lemma~\ref{L2},
$M_\delta$ is dense. It follows that
\begin{align*}
\{ f \in C(\Omega) : \Sigma_{{\rm ac}}(f) = \emptyset \} & = \{ f \in C(\Omega) : M(f) =
0 \} = \bigcap_{\delta > 0} M_\delta
\end{align*}
is residual.
\end{proof}

\begin{proof}[Proof of Theorem~\ref{main2}.]
For $\Lambda, \delta>0$, we define
$$
M_\delta(\Lambda) = \Bigl\{ f \in C(\Omega) : \int_0^\Lambda M(\lambda f) \, d\lambda <
\delta\Bigr\}.
$$
By Lemma~\ref {L1}, $M_\delta(\Lambda)$ is open and by Lemma~\ref {L2},
$M_\delta(\Lambda)$ is dense. Thus,
$$
\bigcap_{\Lambda,\delta > 0} M_\delta(\Lambda)
$$
is residual.  It follows that for Baire generic $f \in C(\Omega)$, we have $\Sigma_{{\rm
ac}}(\lambda f) = \emptyset$ for almost every $\lambda > 0$.
\end{proof}

\section{Concluding Remarks}

\begin{rem}

It is possible to improve Lemma~\ref {L1} to show that $M(f)$ is an upper semi-continuous
function of $f \in L^1(\Omega)$.  The additional point is that, given $f \in L^1(\Omega)$
and $\varepsilon>0$, we can choose a bounded interval $I \subset \R$ such that for every
$\tilde f \in L^1(\Omega)$ close to $f$, we have $|\{E \in \R \setminus I :
\gamma_{\tilde f}(E)=0\}|<\varepsilon$. To see this, one shows first that the integrated
density of states $N_f(E) \in L^\infty(\R)$ is a continuous function of $f \in
L^1(\Omega)$, and then one uses \cite {DS} to bound the size of the absolutely continuous
spectrum near infinity.

\end{rem}

\begin{rem}

By the Wonderland theorem \cite{s2} (see also \cite{ls}), the set of $f$'s leading to
purely singular spectrum is a $G_\delta$ set in all metric topologies that imply strong
resolvent convergence of the associated operators. This permits one to deduce generic
singular spectrum if one can exhibit a dense set with this property. With the Kotani
result (combined with Lemma~\ref{L3}) as input, this only gives a generic set in
$L^\infty(\Omega)$ and does not imply Theorem~\ref{main}. It is not clear how to prove
Theorem~\ref{main} using this strategy, but it would be interesting to find an explicit
dense set of continuous functions such that the corresponding operators have empty
absolutely continuous spectrum.

\end{rem}
\begin{rem}

The result of this paper naturally extends to the context of more general $\SL(2,\R)$
cocycles.  A possible formulation is the following.  Given $A \in C(\Omega,\SL(2,\R))$,
one can consider a one-parameter family of cocycles $(T,R_\theta A)$, where $R_\theta=\bm
\cos \theta&-\sin \theta\\\sin \theta&\cos \theta \em$, and the result is that for
generic $A$ and for almost every $\theta \in \R$, the Lyapunov exponent of $(T,R_\theta
A)$ is positive.  The key point is that the relevant part of Kotani's Theory (which is
used in the proof of Lemma \ref {L2}) can be carried out in this setting (see \cite {AK}
for related results).  (To prove the analogue of Lemma \ref {L1}, one can use \cite {AB}
to show that the average Lyapunov exponent of the family $\theta \mapsto (T,R_\theta A)$
depends continuously on $A$ in the $L^1$ topology.)

Notice that for certain choices of $T$ (say, irrational rotations), there are open sets
$U \subset C(\Omega,\SL(2,\R))$ such that the Lyapunov exponent of $(T,R_\theta A)$ is $0$ for
generic $A \in U$ and for generic
$\theta \in \R$, \cite {B}.  Based on this, some authors have argued that under
weak smoothness requirements (such as continuity), positive Lyapunov exponents are rare.
Our result shows in a sense that positive Lyapunov exponents tend to prevail in a mixed
topological/measure-theoretic category, even when they are topologically rare.

\end{rem}


\begin{thebibliography}{10}

\bibitem{AB} A.\ Avila and J.\ Bochi, A formula with some applications to the
theory of Lyapunov exponents, \textit {Israel J.\ Math.} {\bf 131} (2002), 125--137

\bibitem{AK} A.\ Avila and R.\ Krikorian,  Quasiperiodic $\SL(2,\R)$ cocycles,
In preparation

\bibitem{B} J.\ Bochi, Genericity of zero Lyapunov exponents, \textit{Ergodic Theory
Dynam.\ Systems} {\bf 22} (2002), 1667--1696

\bibitem{bj} J.\ Bourgain and S.\ Jitomirskaya, Absolutely continuous spectrum for 1D quasiperiodic
operators, \textit{Invent.\ Math.} {\bf 148} (2002), 453--463

\bibitem{cfks} H.\ L.\ Cycon, R.\ G.\ Froese, W.\ Kirsch, and B.\ Simon, \textit{Schr\"odinger
Operators with Application to Quantum Mechanics and Global Geometry}, Texts and
Monographs in Physics, Springer-Verlag, Berlin (1987)

\bibitem{dk} D.\ Damanik and R.\ Killip, Ergodic potentials with a discontinuous sampling
function are non-deterministic, Preprint (2004)

\bibitem{DS} P.\ Deift and B.\ Simon, Almost periodic Schr\"odinger operators. III.~The absolutely
continuous spectrum in one dimension, \textit{Commun.\ Math.\ Phys.} {\bf 90} (1983),
389--411

\bibitem{dt} T.\ A.\ Driscoll and L.\ N.\ Trefethen, \textit{Schwarz-Christoffel Mapping}, Cambridge Monographs
on Applied and Computational Mathematics {\bf 8}, Cambridge
University Press, Cambridge (2002)

\bibitem{i} K.\ Ishii, Localization of eigenstates and transport phenomena in one-dimensional
disordered systems, \textit{Suppl.\ Prog.\ Theor.\ Phys.} {\bf 53} (1973), 77--138

\bibitem{k} S.\ Kotani, Ljapunov indices determine absolutely continuous spectra of stationary random
one-dimensional Schr\"odinger operators, in \textit{Stochastic
analysis} (Katata/Kyoto, 1982), pp.~225--247, North-Holland Math.\
Library {\bf 32}, North-Holland, Amsterdam (1984)

\bibitem{k1} S.\ Kotani, Jacobi matrices with random potentials taking finitely many values,
\textit{Rev.\ Math.\ Phys.} {\bf 1} (1989), 129--133

\bibitem{ks} H.\ Kunz and B.\ Souillard, Sur le spectre des op\'erateurs aux diff\'erences finies
al\'eatoires, \textit{Commun.\ Math.\ Phys.} {\bf 78} (1980),
201--246

\bibitem{ls} D.\ Lenz and P.\ Stollmann, Generic sets in spaces of measures and generic singular
continuous spectrum for Delone Hamiltonians, Preprint (2004)

\bibitem{p} L.\ Pastur, Spectral properties of disordered systems in one-body approximation,
\textit{Commun.\ Math.\ Phys.} {\bf 75} (1980), 179--196

\bibitem{s} B.\ Simon, Kotani theory for one-dimensional stochastic Jacobi matrices, \textit{Commun.\
Math.\ Phys.} {\bf 89} (1983), 227--234

\bibitem{s2} B.\ Simon, Operators with singular continuous spectrum. I.~General operators, \textit{Ann.\
of Math.} {\bf 141} (1995), 131--145

\end{thebibliography}
\end{document}